\documentclass[11pt]{article}
\usepackage{latexsym,bm}
\usepackage{mathrsfs}
\usepackage{amsmath,amsthm}
\usepackage{graphicx}
\usepackage{tikz}
\usepackage{needspace}
\usepackage{amssymb}
\usepackage{color}
\usepackage[noadjust]{cite}
\usepackage{array}
\usepackage{stfloats}
\usepackage[colorlinks,
            linkcolor=blue,
            anchorcolor=blue,
            citecolor=blue,
            ]{hyperref}
\usepackage{caption}
\usepackage{upgreek}

\theoremstyle{plain}
\newtheorem{thm}{Theorem}[section]

\newtheorem{lem}[thm]{Lemma}
\newtheorem{prop}[thm]{Proposition}

\theoremstyle{definition}
\newtheorem{defn}[thm]{Definition}

\theoremstyle{plain}

\theoremstyle{definition}
\newtheorem{prob}{Problem}

\theoremstyle{plain}

\theoremstyle{plain}

\theoremstyle{plain}

\theoremstyle{plain}

\usepackage[top=2.5cm,bottom=2.5cm,left=2.8cm,right=2.2cm]{geometry}
\normalsize \rm
\allowdisplaybreaks[4]
\makeatletter
\@addtoreset{equation}{section}
\makeatother

 \usepackage{indentfirst}
\usepackage{enumitem}
\setlist[enumerate]{leftmargin=2.7em,itemsep=0.35em,topsep=0.4em}
\hypersetup{pdftitle={The Hamilton cycle space in the Erd\H{o}s--R\'enyi random graph process}}
\newcommand{\Ftwo}{\mathbb F_2}
\newcommand{\CS}{\mathcal C}
\newcommand{\HS}{\mathcal H}

\newcommand{\Prb}{\mathbb P}
\newcommand{\Ex}{\mathbb E}
\newcommand{\Var}{\operatorname{Var}}
\newcommand{\sstar}{s_*}
\newcommand{\bstar}{b_*}
\newcommand{\hstar}{h_*}
\newcommand{\qstar}{q_*}

\begin{document}
\renewcommand{\figurename}{Fig.}
\title{The Hamilton cycle space in the
Erd\H{o}s--R\'enyi random graph process}
\author{
Hao Zhang$^{1}$\quad  Wensheng Sun$^{2}$\quad  Panpan Cheng$^{1}$%
\thanks{Corresponding author.}\\[0.4em]
\begin{minipage}{0.92\textwidth}
\centering
\small
$^{1}$ School of Mathematics, Shandong University, Jinan 250100, China\\
$^{2}$ School of Mathematics and Statistics, Gansu Center for Applied Mathematics,\\
Lanzhou University, Lanzhou, Gansu 730000, China
\end{minipage}
}

\date{}

\footnotetext[0]{E-mail addresses: zhang.h@mail.sdu.edu.cn (H. Zhang), wensheng07002@163.com (W. Sun), ppancheng@126.com (P. Cheng).}
\maketitle
\vspace{1mm}

\begin{abstract}
Let $G(n,m)$ denote the graph at step $m$ of the Erd\H{o}s--R\'enyi random graph process on $n$ vertices. The process starts with the edgeless graph and adds one edge at each step, chosen uniformly at random from the missing edges. Let $\tau_3$ be the first step at which the minimum degree of $G(n,m)$ reaches three. We answer a question of Christoph, Nenadov and Petrova [J. Combin. Theory Ser. B 176 (2026) 254--267] by proving that, as $n$ tends to infinity through odd integers, with high probability the Hamilton cycles of $G(n,m)$ span its cycle space over $\Ftwo$ simultaneously for all $m\ge\tau_3$. On the same event, $G(n,m)$ is Hamilton-connected for every $m\ge\tau_3$.


\end{abstract}

\section{Introduction}\label{sec1}

The systematic study of random graphs began with the work of Erd\H{o}s and R\'enyi~\cite{erd1,erd2}. They considered a graph chosen uniformly at random from all graphs with $m$ edges on a fixed set of $n$ labeled vertices and studied how properties such as connectivity, planarity and Hamiltonicity depend on $m$. A classical result of Erd\H{o}s and R\'enyi~\cite{erd1} gives the limiting probability of connectivity near the threshold. More precisely, if $m=\lfloor n(\log n+c)/2\rfloor$ for a fixed real number $c$, then the probability that the graph is connected tends to $\exp(-e^{-c})$ as $n\to\infty$. Another classical result concerns the first appearance of a Hamilton cycle in the random graph process. Ajtai, Koml\'os and Szemer\'edi~\cite{ajt} and, independently, Bollob\'as~\cite{bol1} proved that, with high probability, the first Hamilton cycle appears precisely when the minimum degree reaches two. Random graph theory has since developed into a major branch of probabilistic combinatorics. For further results on random graphs, particularly on Hamilton cycles, we refer the reader to~\cite{gle,mon,chr}.

Throughout the paper, let $G(n,m)$ denote the graph at step $m$ of the Erd\H{o}s--R\'enyi random graph process on $[n]=\{1,\ldots,n\}$. The initial graph $G(n,0)$ has no edges. For $0\le m<\binom{n}{2}$, we obtain $G(n,m+1)$ from $G(n,m)$ by adding an edge chosen uniformly at random from the missing edges. Equivalently, we may expose the edges of $K_n$ in a uniformly random order and let $G(n,m)$ consist of the first $m$ edges.  Thus, for each fixed $m$, the graph $G(n,m)$ is uniformly distributed over all graphs on $[n]$ with $m$ edges, and this process ends with $G(n,\binom{n}{2})=K_n$. We also use the binomial model of Gilbert~\cite{gil}, whose distribution we denote by $\mathbb{G}(n,p)$. In this model, the vertex set is $[n]$ and each possible edge is present independently with probability $p$. For $0<p<1$, a graph sampled from this model, conditional on having exactly $m$ edges, has the same distribution as $G(n,m)$. Since the same random ordering of the edges determines $G(n,m)$ for every $m$, we can study when a property first appears in the process and whether it persists at all subsequent steps.

We say that an event holds \emph{with high probability} if its probability tends to one as $n\to\infty$. For $1\le k\le n-1$, let $\tau_k=\min\{m:\delta(G(n,m))\ge k\}$, where $\delta(G)$ denotes the minimum degree of $G$. For every fixed positive integer $k$, Bollob\'as and Thomason~\cite{bol2} proved that, with high probability, the graph first becomes $k$-connected at $\tau_k$. A related result of Ajtai, Koml\'os and Szemer\'edi~\cite{ajt} and, independently, Bollob\'as~\cite{bol1} states that the first Hamilton cycle appears at $\tau_2$ with high probability. Thus, with high probability, each property first appears as soon as the corresponding minimum-degree condition is met. Since both $k$-connectivity and Hamiltonicity are monotone increasing under edge addition, each persists at every subsequent step.

We study the subspace of the cycle space generated by Hamilton cycles.  For a graph $G$, we identify each edge set with its incidence vector in $\Ftwo^{E(G)}$. The \emph{cycle space} $\CS(G)$ is spanned by the cycles of $G$, and the \emph{Hamilton space} $\HS(G)$ is the subspace spanned by its Hamilton cycles. We say that $G$ is \emph{Hamilton-generated} if $\HS(G)=\CS(G)$. Equivalently, every cycle of $G$ can be expressed as a symmetric difference of Hamilton cycles. In 1990, Alspach, Locke and Witte~\cite{als} determined the Hamilton space of every connected Cayley graph on a finite abelian group. In particular, every such graph of odd order is Hamilton-generated. For random graphs, a related question concerns the subspace spanned by cycles of a fixed odd length. DeMarco, Hamm and Kahn~\cite{dem} proved that, with high probability, the triangles of a binomial random graph span its cycle space whenever every edge lies in a triangle.  Subsequently, Baron and Kahn~\cite{bar}  extended this result to every fixed odd length $k\ge3$. For Hamilton cycles, the relevant length is $n$ itself and thus grows with the order of the graph.

In particular, the parity of the order places a restriction on Hamilton generation. If $G$ has even order, every symmetric difference of Hamilton cycles has even size, so an odd cycle prevents Hamilton generation. Now suppose that $G$ has odd order and contains a vertex $v$ of degree two. If $G-v$ contains an odd cycle $C$, then $C$ cannot be expressed as a symmetric difference of Hamilton cycles of $G$: every Hamilton cycle contains both edges incident with $v$, while $C$ contains neither. Thus, any such representation would require an even number of Hamilton cycles for these edges to cancel. Since each Hamilton cycle has odd length, the resulting symmetric difference would have an even number of edges, contrary to the choice of $C$. Heinig~\cite{hei2} used this observation to relate Hamilton generation in the binomial random graph to the threshold for minimum degree three.

For odd $n$, Heinig~\cite{hei2} proved that a graph sampled from $\mathbb{G}(n,p)$ is Hamilton-generated with high probability if $p\ge n^{-1/2+\varepsilon}$ for any fixed $\varepsilon>0$. Recently, Christoph, Nenadov and Petrova~\cite{chr} improved this bound to $p\ge C\log n/n$ for a sufficiently large constant $C$, as a consequence of a theorem for pseudorandom graphs. Hefetz and Krivelevich~\cite{hef} then obtained the result for $p\ge(\log n+2\log\log n+\omega(1))/n$, at the threshold for minimum degree three in the binomial model. This determines the minimum degree threshold in the binomial model. However, the corresponding hitting-time result for the random graph process $G(n,m)$ remained open. This led Christoph, Nenadov and Petrova~\cite{chr} to pose the following question.

\begin{prob}[\cite{chr}]\label{pro1}
Suppose that $n$ tends to infinity through odd integers. Is it true that, with high probability, $\HS(G(n,m))=\CS(G(n,m))$ for every $m\ge\tau_3$?
\end{prob}

A graph is \emph{Hamilton-connected} if any two distinct vertices are joined by a spanning path. By Heinig's lemma~\cite{hei1}, if a graph is both Hamilton-connected and Hamilton-generated, then it retains both properties under edge addition. It therefore suffices to prove that, with high probability, $G(n,\tau_3)$ is both Hamilton-connected and Hamilton-generated. We  give the main result as follows.

\begin{thm}\label{the1.1}
As $n$ tends to infinity through odd integers, with high probability $G(n,m)$ is Hamilton-connected and satisfies $\HS(G(n,m))=\CS(G(n,m))$ simultaneously for all $m\ge\tau_3$.
\end{thm}

We use two standard probabilistic tools; see Alon and Spencer~\cite{alo}. The \emph{Chernoff bounds} control deviations of a sum $X$ of independent Bernoulli random variables from its mean $\mu=\Ex X>0$. In particular, for $0<\varepsilon<1$,
\begin{equation}\label{eq-chernoff}
\begin{aligned}
\Prb\bigl(X\le(1-\varepsilon)\mu\bigr)&\le\exp(-\varepsilon^2\mu/2)\text{ and }
\Prb\bigl(|X-\mu|\ge\varepsilon\mu\bigr)\le2\exp(-\varepsilon^2\mu/3).
\end{aligned}
\end{equation}
We also use the one-sided estimates
\begin{equation}\label{eq-chernoff-one-sided}
\begin{aligned}
\Prb(X\le t)&\le\exp\bigl(-\mu+t+t\log(\mu/t)\bigr)
&& (0<t<\mu),\\
\Prb(X\ge t)&\le(e\mu/t)^t
&& (t\ge\mu).
\end{aligned}
\end{equation}
The first estimate bounds the probability that $X$ is substantially smaller than its mean, and the second bounds the probability that it is larger.

We next recall the classical \textit{Lov\'asz local lemma} (LLL) in its symmetric form~\cite{alo}. It guarantees that a family of events can be avoided simultaneously if each event has small probability and depends on sufficiently few of the others. More precisely, let $E_1,\ldots,E_t$ be events, each of probability at most $q$. Suppose that, for each $i$, one can omit at most $D$ of the other events so that $E_i$ is independent of the remaining events taken together. If $eq(D+1)\le 1$, then
\[
\Prb\left(\bigcap_{i=1}^t E_i^{\mathrm c}\right)>0.
\]
We use this lemma to partition a vertex set $Y$ into two parts so that, for each vertex under consideration, each part contains at least a fixed positive proportion of its neighbors in $Y$.

Section~\ref{sec2} gives a deterministic \emph{criterion}, consisting of conditions~\ref{iteA1}--\ref{iteA5} in Definition~\ref{def2.1}, that allows for vertices of small degree. We prove that these conditions imply Hamilton-connectedness in Section~\ref{sec3} and $\HS(G)=\CS(G)$ for odd-order graphs in Section~\ref{sec4}. In Section~\ref{sec5}, we verify the criterion at $\tau_3$ and show that both conclusions persist under further edge additions.

\section{The deterministic criterion and auxiliary results}\label{sec2}

Vertices of small degree need separate treatment, since they need not satisfy the expansion estimates used below. We allow an exceptional set with sufficiently separated neighborhoods. For the remaining vertices, lower degree bounds and upper bounds on edge counts in small sets imply the required expansion.

For a graph $G$ and a set $X\subseteq V(G)$, let
\[
N_G(X)=\{w\in V(G)\setminus X:vw\in E(G)\text{ for some }v\in X\}
\]
be the external neighborhood of $X$, and let $N_G(v)=N_G(\{v\})$. Write $d_G(v)=|N_G(v)|$ and $d_G(v,X)=|N_G(v)\cap X|$. Let $e_G(X)$ denote the number of edges with both ends in $X$ and, for disjoint sets $X,Y\subseteq V(G)$, let $e_G(X,Y)$ denote the number of edges between $X$ and $Y$. The length of a path is its number of edges. For $v\in V(G)$ and a nonnegative integer $h$, the \emph{ball} of radius $h$ around $v$ is
\[
B_G(v,h)=\{w\in V(G):\operatorname{dist}_G(v,w)\le h\},
\]
where $\operatorname{dist}_G(v,w)$ is the length of a shortest $v,w$-path and is infinite if no such path exists. All logarithms are natural. We assume throughout that $n$ is sufficiently large in terms of the fixed constants in each statement.

We formulate the criterion using the typical properties of sparse random graphs identified by Hefetz and Krivelevich~\cite{hef}. Let $L=\log n$, $\ell=\log\log n$, $\qstar=L/\ell$ and $r=\lceil\sqrt L\rceil$, and let
\begin{equation}\label{eq-1}
\begin{gathered}
\sstar=\left\lfloor\frac{n\ell^2}{L}\right\rfloor,
\; \bstar=\left\lceil\frac{n\ell^{3/2}}{L}\right\rceil,\;
p_0=\frac{L+2\ell}{n}\text{ and } \hstar=\left\lceil\frac{6L}{\ell}\right\rceil.
\end{gathered}
\end{equation}
The parameter $\qstar$ controls the edge counts in small sets, $\hstar$ bounds the lengths of the paths used below, and $p_0$ is a reference edge density; probabilities of events are denoted by $\Prb$. The parameters $L,\ell,\qstar,r$ and those in~\eqref{eq-1} depend on the order $n$ of the ambient graph, including when they occur in statements about subgraphs. We shall use the relations $\bstar=o(\sstar)$, $r\sstar=o(n)$ and $\hstar=o(L)$.

\begin{defn}\label{def2.1}
Let $G$ be an $n$-vertex graph, and let $S\subseteq V(G)$. Let $T=N_G(S)$, $Z=S\cup T$ and $U=V(G)\setminus Z$. We say that $(G,S)$ satisfies the \emph{criterion} if conditions~\ref{iteA1}--\ref{iteA5} below hold.
\begin{enumerate}[label=(A\arabic*),ref=(A\arabic*)]
\item\label{iteA1} $\delta(G)\ge3$, $\Delta(G)\le10L$, and $d_G(v)\ge L/10$ for every $v\notin S$.
\item\label{iteA2} $|Z|\le\sqrt n$, both $S$ and $T$ are independent, and $d_G(v,Z)\le1$ for every $v\notin S$.
\item\label{iteA3} $e_G(X)\le\qstar|X|$ whenever $|X|\le\sstar$.
\item\label{iteA4} $e_G(X,Y)\le\qstar|X|$ whenever $X,Y$ are disjoint, $|X|\le\sstar$, and $|Y|\le r|X|$.
\item\label{iteA5} For any two disjoint vertex sets $X,Y$ with $|X|,|Y|\ge\bstar$, we have
\begin{equation}\label{eq-2}
0.999p_0|X||Y|\le e_G(X,Y)\le1.001p_0|X||Y|.
\end{equation}
\end{enumerate}
\end{defn}

The lower degree bound in~\ref{iteA1} is three on $S$ and $L/10$ outside $S$. In Section~\ref{sec5}, we choose $S$ from the graph at the lower endpoint of an interval of the random graph process. The same choice then works for every graph in that interval.

\begin{thm}\label{the2.2}
Let $G$ be an $n$-vertex graph and $S\subseteq V(G)$. If $(G,S)$ satisfies conditions~\ref{iteA1}--\ref{iteA5}, then $G$ is Hamilton-connected and, if $n$ is odd, satisfies $\HS(G)=\CS(G)$.
\end{thm}

For Theorem \ref{the2.2}, we prove Hamilton-connectedness in Section~\ref{sec3} and then use it in the proof of Hamilton generation in Section~\ref{sec4}. We first recall two results on expansion and disjoint paths.

\begin{defn}\label{def2.3}
A graph $G$ on $n\ge3$ vertices is a \emph{$c$-expander} if $|N_G(X)|\ge c|X|$ for every $X\subseteq V(G)$ with $|X|<n/(2c)$ and every two disjoint vertex sets of size at least $n/(2c)$ have an edge between them.
\end{defn}

\begin{thm}[{\cite{dra2}}]\label{the2.4}
There is a constant $c_0>0$ such that every $c$-expander with $c\ge c_0$ is Hamilton-connected.
\end{thm}

Let $m_0,d_0$ be positive integers, and let $\alpha\in(0,1)$. A graph $G$ has \emph{property $P_\alpha(m_0,d_0)$} if $|N_{G-\mathcal F}(X)|\ge2d_0|X|$ for every $X\subseteq V(G)$ with $|X|\le m_0$ and every $\mathcal F\subseteq E(G)$ such that $d_{\mathcal F}(v)\le\alpha d_G(v)$ for each $v\in X$. Here $d_{\mathcal F}(v)$ is the number of edges of $\mathcal F$ incident with $v$.

The next result is a consequence of the routing theorem of Dragani\'c, Krivelevich and Nenadov~\cite{dra1}. We use the form stated by Hefetz and Krivelevich~\cite{hef}.

\begin{thm}[{\cite{hef}}]\label{the2.5}
Let $m_0,d_0$ be integers with $3\le d_0<m_0$, and let $1/2<\alpha<1$ and $0<\beta<2\alpha-1$. Suppose that $G$ has property $P_\alpha(m_0,d_0)$ and that every two disjoint vertex sets of size at least $m_0(d_0-1)/16$ have an edge between them. Let $\mathcal T\subseteq V(G)$ satisfy $d_G(v,\mathcal T)\le\beta d_G(v)$ for every $v\in V(G)$. Then any collection of at most $d_0m_0\log d_0/(15\log m_0)$ pairwise disjoint pairs of vertices in $\mathcal T$ can be joined by pairwise vertex-disjoint paths in $G$.
\end{thm}

\begin{lem}\label{lem2.6}
Let $G$ be an $n$-vertex graph satisfying~\ref{iteA3} and~\ref{iteA4}, and let $\gamma>0$ and $\alpha\in[0,1)$ be fixed. If $H\subseteq G$ has minimum degree at least $\gamma L$, then $|N_{H-\mathcal F}(X)|\ge r|X|$ for every $X\subseteq V(H)$ with $|X|\le\sstar$ and every $\mathcal F\subseteq E(H)$ such that $d_{\mathcal F}(v)\le\alpha d_H(v)$ for each $v\in X$.
\end{lem}
\begin{proof}
We argue by contradiction. Suppose that there are a nonempty set $X\subseteq V(H)$ with $|X|\le\sstar$ and an edge set $\mathcal F\subseteq E(H)$ satisfying the stated degree condition such that $|N_{H-\mathcal F}(X)|<r|X|$. Each vertex of $X$ retains at least a proportion $1-\alpha$ of its incident edges after the deletion of $\mathcal F$. On the other hand, \ref{iteA3} and \ref{iteA4} bound the edges within $X$ and between $X$ and its external neighborhood. Thus
\[
(1-\alpha)\gamma L|X|
\le\sum_{v\in X}d_{H-\mathcal F}(v)
\le2e_G(X)+e_G\bigl(X,N_{H-\mathcal F}(X)\bigr)
\le3\qstar|X|.
\]
This contradicts $\qstar/L=1/\ell\to0$ for sufficiently large $n$.
\end{proof}

\begin{lem}\label{lem2.7}
Let $G$ be an $n$-vertex graph satisfying~\ref{iteA3}--\ref{iteA5}, and let $\eta,\gamma>0$ be fixed. If $W\subseteq V(G)$ satisfies $|W|\ge\eta n$ and $\delta(G[W])\ge\gamma L$, then $G[W]$ is Hamilton-connected.
\end{lem}
\begin{proof}
Choose a fixed $c\ge\max\{c_0,3\}$. For large $n$ we have $r\ge c$, so by Lemma~\ref{lem2.6}, the expansion bound holds for sets of size at most $\sstar$. Now suppose that $\sstar<|X|<|W|/(2c)$ and $|N_{G[W]}(X)|<c|X|$. Then $W\setminus\bigl(X\cup N_{G[W]}(X)\bigr)$ has size greater than $|W|(c-1)/(2c)\ge |W|/3$ and sends no edges to $X$. Both sets have size at least $\bstar$, contrary to~\ref{iteA5}. Moreover, any two disjoint subsets of $W$ of size at least $|W|/(2c)$ have an edge between them by~\ref{iteA5}. The graph $G[W]$ is therefore a $c$-expander, so we may apply Theorem~\ref{the2.4}.
\end{proof}

\begin{lem}\label{lem2.8}
Let $G$ be an $n$-vertex graph satisfying~\ref{iteA3}--\ref{iteA5}, and let $\gamma>0$ be fixed. Let $H$ be obtained from an induced subgraph of $G$ by deleting an edge set of maximum degree at most four, and suppose that $\delta(H)\ge\gamma L$. If $\mathcal T\subseteq V(H)$ satisfies $d_H(v,\mathcal T)\le d_H(v)/10$ for every $v\in V(H)$, then any collection of at most $\sqrt n$ pairwise disjoint pairs of vertices in $\mathcal T$ can be joined by pairwise vertex-disjoint paths in $H$.
\end{lem}
\begin{proof}
We use Theorem~\ref{the2.5} with $m_0=16\lfloor n/(16L)\rfloor$, $d_0=\lfloor r/4\rfloor$, $\alpha=2/3$ and $\beta=1/10$. The inequalities $m_0\le\sstar$ and $2d_0\le r$ allow us to apply Lemma~\ref{lem2.6}, so $H$ has property $P_\alpha(m_0,d_0)$. To check the edge condition, let $h_0=m_0(d_0-1)/16$. We have $h_0\asymp n/\sqrt L\gg\bstar$. If $X,Y\subseteq V(H)$ are disjoint and each has size at least $h_0$, the deletion removes at most $4\min\{|X|,|Y|\}$ edges between them. By~\ref{iteA5},
\[
e_H(X,Y)\ge
\min\{|X|,|Y|\}\bigl(0.999p_0\max\{|X|,|Y|\}-4\bigr)>0,
\]
where the last inequality follows from $p_0 h_0\to\infty$. We also have $\beta<2\alpha-1$. Since the number of pairs allowed by Theorem~\ref{the2.5} is of order $n\ell/L^{3/2}$ and this exceeds $\sqrt n$ for sufficiently large $n$, the theorem applies to every collection in the statement.
\end{proof}

\section{Hamilton-connectedness under (A1)--(A5)}
\label{sec3}
Throughout this section, $(G,S)$ satisfies conditions~\ref{iteA1}--\ref{iteA5} of Definition~\ref{def2.1}. We retain the notation $T=N_G(S)$, $Z=S\cup T$ and $U=V(G)\setminus Z$, and call the vertices outside $S$ \emph{ordinary}. We shall first join paths through the exceptional vertices by disjoint paths in one part of the ordinary vertices. A second part is reserved so that the vertices left over still have sufficiently large degree. By the following partition lemma, we can choose both parts so that each retains enough neighbors of every ordinary vertex.

\begin{lem}\label{lem3.1}
Let $G$ be an $n$-vertex graph with $\Delta(G)\le10L$. Let $\gamma>0$ be fixed, and let $Q,Y\subseteq V(G)$ satisfy $d_G(v,Y)\ge\gamma L$ for every $v\in Q$. Then $Y$ has a partition $Y=Y_1\mathbin{\dot\cup}Y_2$ such that $d_G(v,Y_i)\ge d_G(v,Y)/3\ge\gamma L/3$ for every $v\in Q$ and $i\in\{1,2\}$. Moreover, if $|Q|\ge\eta n$ for some fixed $\eta>0$, then $|Y_i|\ge\eta\gamma n/30$ for each $i\in\{1,2\}$.
\end{lem}

\begin{proof}
Color each vertex of $Y$ independently with either color, each with probability $1/2$. For a vertex $v\in Q$, the Chernoff bound~\eqref{eq-chernoff} shows that the probability that one of the colors appears fewer than $d_G(v,Y)/3$ times in $N_G(v)\cap Y$ is at most $2\exp(-\gamma L/36)$. This event depends only on the colors in $N_G(v)\cap Y$. Each vertex there has at most $\Delta(G)$ neighbors in $Q$, so the event is independent of all but at most $\Delta(G)^2\le100L^2$ of the other events. For sufficiently large $n$, we have $2e\exp(-\gamma L/36)(100L^2+1)<1$. By the symmetric Lov\'asz local lemma stated in Section~\ref{sec1}, we can therefore choose a coloring in which each vertex of $Q$ has at least a third of its neighbors in $Y$ in each color class.

For each color class $Y_i$, the sum $\sum_{v\in Q}d_G(v,Y_i)$ is at least $|Q|\gamma L/3$, while every vertex in $Y_i$ contributes at most $\Delta(G)$ incidences. Thus this number is at most $\Delta(G)|Y_i|\le10L|Y_i|$, and $|Q|\ge\eta n$ implies $|Y_i|\ge\eta\gamma n/30$.
\end{proof}

The following spanning-path statement adapts an argument of Hefetz and Krivelevich~\cite{hef} to prescribed ordinary endpoints. We include the proof. The condition that each exceptional vertex has two neighbors outside the endpoints will also be used in Section~\ref{sec4} to complete two paths to Hamilton cycles with one common complementary path.

\begin{lem}\label{lem3.2}
Let $\eta,\gamma>0$ be fixed, let $W\subseteq V(G)$ satisfy $|W|\ge\eta n$, and let $x,y$ be distinct vertices of $W\setminus S$. Suppose that each $v\in W\setminus S$ has at least $\gamma L$ neighbors in $W$ and that each $u\in W\cap S$ has at least two neighbors in $W\setminus\{x,y\}$. Then $G[W]$ contains a spanning $x,y$-path.
\end{lem}

\begin{proof}
If $W\cap S=\varnothing$, we may apply Lemma~\ref{lem2.7} directly. We may therefore assume that $W\cap S=\{u_1,\ldots,u_t\}$ with $t\ge1$. For each $i$, choose distinct vertices $a_i,b_i\in N_G(u_i)\cap(W\setminus\{x,y\})$. Since $S$ is independent, all these neighbors are ordinary. Moreover, an ordinary vertex has at most one neighbor in $S$ by~\ref{iteA2}, so all $2t$ chosen vertices are distinct.

Let $Z_0=\{u_i,a_i,b_i:1\le i\le t\}$ and $Y=W\setminus(Z_0\cup\{x,y\})$. By~\ref{iteA2}, every ordinary vertex has at most one neighbor in $Z_0\subseteq Z$. Together with the two endpoints, the deletion therefore removes at most three of its neighbors in $W$. Consequently, $d_G(v,Y)\ge\gamma L-3\ge3\gamma L/4$ for every $v\in W\setminus S$ when $n$ is sufficiently large. We may apply Lemma~\ref{lem3.1} to obtain a partition $Y=A\mathbin{\dot\cup}B$ in which every ordinary vertex of $W$ has at least $\gamma L/4$ neighbors in each part. Since $|W\setminus S|\ge\eta n-\sqrt n$, the size bound in that lemma shows that both $A$ and $B$ have order at least a positive constant times $n$.

We use $A$ to join the paths through the exceptional vertices, leaving $B$ untouched. Let $\mathcal T=\{b_1,\ldots,b_t,a_2,\ldots,a_t,y\}$ and $H=G[A\cup\mathcal T]$. Every vertex of $H$ is ordinary and has at least $\gamma L/4$ neighbors in $A$. Since the chosen $a_i,b_i$ all lie in $Z$, each vertex of $H$ has at most one neighbor among them by~\ref{iteA2}. The only other possible neighbor in $\mathcal T$ is $y$. Thus $d_H(v,\mathcal T)\le2\le d_H(v)/10$ for every $v\in V(H)$ and sufficiently large $n$. There are $t\le|S|\le\sqrt n$ pairs to join, so by Lemma~\ref{lem2.8}, we can find pairwise vertex-disjoint paths joining $b_t$ to $y$ and $b_{i-1}$ to $a_i$ for each $2\le i\le t$.

Let $W'$ consist of the vertices of $W$ left after deleting $u_1,\ldots,u_t$ and all vertices on these paths. In particular, $x$, $a_1$ and all of $B$ remain in $W'$. Since every vertex of $W'$ is ordinary and has at least $\gamma L/4$ neighbors in $B$, we may apply Lemma~\ref{lem2.7} to $G[W']$. Choose a spanning $x,a_1$-path in this graph and continue it along $a_1u_1b_1$. At $b_i$, for $i=1,\ldots,t-1$, follow the chosen path to $a_{i+1}$ and then traverse $a_{i+1}u_{i+1}b_{i+1}$. Finally, follow the chosen $b_t,y$-path. The disjointness of the connecting paths and the choice of $W'$ ensure that this is a path and that it visits every vertex of $W$. Figure~\ref{fig2} shows the construction for $t=2$. The dashed paths have their internal vertices in $A$, while the $x,a_1$-path spans $W'$, including all of $B$ and the unused vertices of $A$.

\par\medskip
\noindent\begin{minipage}{\linewidth}
\centering
\resizebox{0.95\linewidth}{!}{%
\begin{tikzpicture}[
  ordinary/.style={circle,draw=black,fill=white,minimum size=6.5mm,inner sep=0pt,font=\small},
  exceptional/.style={ordinary,fill=black,text=white},
  point/.style={circle,draw=black,fill=white,inner sep=1.2pt,line width=0.5pt},
  route/.style={line width=0.85pt,densely dashed},
  spanning/.style={line width=1pt},
  every node/.style={font=\small}
]
\draw[rounded corners=3pt,fill=black!3,line width=0.4pt] (-0.25,0.05) rectangle (4.55,3.05);
\node at (2.15,2.74) {$W'$};
\draw[fill=black!12,line width=0.3pt] (1.12,1.5) ellipse (0.73 and 0.99);
\draw[fill=white,line width=0.3pt] (2.85,1.5) ellipse (0.73 and 0.99);
\node at (1.12,0.32) {$B$};
\node at (2.85,0.32) {$A\cap W'$};
\node[ordinary] (x) at (0.15,1.5) {$x$};
\node[ordinary] (a1) at (4.12,1.5) {$a_1$};
\node[exceptional] (u1) at (5.25,1.5) {$u_1$};
\node[ordinary] (b1) at (6.35,1.5) {$b_1$};
\node[ordinary] (a2) at (8.05,1.5) {$a_2$};
\node[exceptional] (u2) at (9.15,1.5) {$u_2$};
\node[ordinary] (b2) at (10.25,1.5) {$b_2$};
\node[ordinary] (y) at (12.05,1.5) {$y$};
\coordinate (p1) at (0.82,1.13);
\coordinate (p2) at (1.43,2.05);
\coordinate (p3) at (2.52,1.13);
\coordinate (p4) at (3.13,2.05);
\draw[spanning] (x) -- (p1) -- (p2) -- (p3) -- (p4) -- (a1);
\foreach \p in {p1,p2,p3,p4} \node[point] at (\p) {};
\draw[spanning] (a1) -- (u1) -- (b1);
\draw[spanning] (a2) -- (u2) -- (b2);
\coordinate (r1) at (6.82,2.28);
\coordinate (r2) at (7.58,2.28);
\coordinate (r3) at (10.75,2.28);
\coordinate (r4) at (11.55,2.28);
\draw[route] (b1) -- (r1) -- (r2) -- (a2);
\draw[route] (b2) -- (r3) -- (r4) -- (y);
\foreach \p in {r1,r2,r3,r4} \node[point] at (\p) {};
\node[font=\footnotesize,align=center] at (7.2,2.85) {internal vertices\ in $A$};
\node[font=\footnotesize,align=center] at (11.15,2.85) {internal vertices\ in $A$};
\node[font=\footnotesize] at (8.1,0.52) {$u_1,u_2\in S$};
\end{tikzpicture}%
}
\captionsetup{font=small,labelsep=period,hypcap=false}

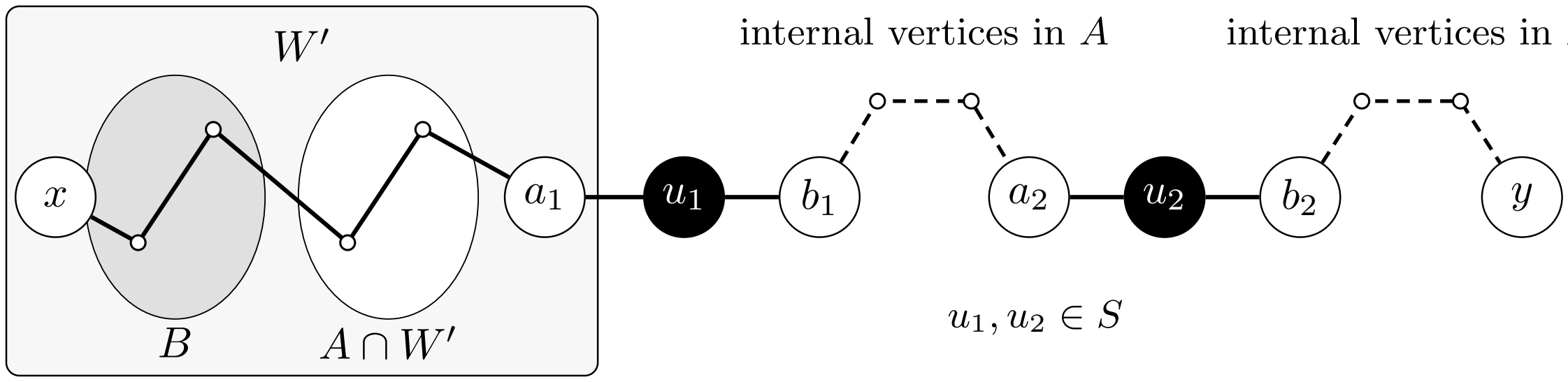
\captionof{figure}{The spanning path in Lemma~\ref{lem3.2} for $t=2$.}
\label{fig2}
\end{minipage}
\par\medskip

\end{proof}

\Needspace{5\baselineskip}
\begin{prop}\label{pro3.3}
The graph $G$ is Hamilton-connected.
\end{prop}
\begin{proof}
Let $x,y$ be distinct vertices of $G$. We shall use Lemma~\ref{lem3.2} with $\eta=1/2$ and $\gamma=1/20$. The sets $W$ considered below omit at most two vertices of $G$. Thus every ordinary vertex in $W$ has at least $L/10-2\ge L/20$ neighbors there, as required by the lemma.

Suppose first that $x,y\notin S$. If each vertex of $S$ has at least two neighbors outside $\{x,y\}$, we apply Lemma~\ref{lem3.2} with $W=V(G)$. Otherwise, the bound $\delta(G)\ge3$ implies that some $u\in S$ has exactly three neighbors $x,y,z$. By~\ref{iteA2}, $z$ is ordinary and $u$ is the only neighbor in $S$ of each of $x,y,z$. Hence every vertex of $S\setminus\{u\}$ retains its neighbors when $x,u,z,y$ are excluded. We may therefore apply Lemma~\ref{lem3.2} with $W=V(G)\setminus\{x,u\}$ and endpoints $z,y$. The resulting spanning $z,y$-path in $G[W]$, preceded by $xuz$, is a spanning $x,y$-path in $G$.

Now suppose that exactly one of $x,y$ belongs to $S$, say $x\in S$, and choose $a\in N_G(x)\setminus\{y\}$. The vertex $a$ is ordinary and has $x$ as its only neighbor in $S$ by~\ref{iteA2}. Since $S$ is independent, the deletion of $x,a,y$ removes at most one neighbor of each vertex in $S\setminus\{x\}$. Each such vertex therefore has at least two neighbors left by~\ref{iteA1}. We may apply Lemma~\ref{lem3.2} in $G-x$ with endpoints $a,y$ and then extend the resulting spanning path by the edge $xa$.

If both $x,y\in S$, choose $a\in N_G(x)$ and $b\in N_G(y)$. Condition~\ref{iteA2} implies that $a,b$ are distinct and ordinary, with no neighbor in $S\setminus\{x,y\}$. As $S$ is independent, every remaining exceptional vertex retains all its neighbors in $V(G)\setminus\{x,y,a,b\}$. Apply Lemma~\ref{lem3.2} to $G-\{x,y\}$ with endpoints $a,b$, and extend the resulting path at its two ends by $xa$ and $by$.
\end{proof}

\section{Hamilton cycles span the cycle space}
\label{sec4}

Throughout this section, $(G,S)$ satisfies conditions~\ref{iteA1}--\ref{iteA5} and $n$ is odd. We retain the notation $T,Z,U$ from Definition~\ref{def2.1} and prove that $\HS(G)=\CS(G)$.

We identify edge sets with their incidence vectors in $\Ftwo^{E(G)}$. For $A,B\subseteq E(G)$, define
\[
\langle A,B\rangle=|A\cap B|\pmod2, \  \mbox{and} \ A\perp B\ \Longleftrightarrow\ \langle A,B\rangle=0.
\]
Thus $A$ and $B$ are \emph{orthogonal} if their intersection has even cardinality. For a subspace $\mathcal W\subseteq\Ftwo^{E(G)}$, its \emph{orthogonal complement} is
\[
\mathcal W^\perp=\{A\in\Ftwo^{E(G)}:\langle A,B\rangle=0
\text{ for every }B\in\mathcal W\}.
\]
An \emph{edge cut} is a set of the form
\[
\partial_G X=\{uv\in E(G):u\in X,\ v\in V(G)\setminus X\}
\qquad(X\subseteq V(G)).
\]
The following lemma of Christoph, Nenadov and Petrova~\cite{chr} provides a subgraph that detects a parity obstruction to $\HS(G)=\CS(G)$. We include its proof.

\begin{lem}\label{lem4.1}
If $\HS(G)\ne\CS(G)$, then there is a spanning subgraph $R\subsetneq G$ such that every Hamilton cycle of $G$ meets $E(R)$ in an even number of edges, $E(R)$ is not an edge cut of $G$, and
\begin{equation}\label{eq-3}
e_R(X,V(G)\setminus X)\ge\frac12e_G(X,V(G)\setminus X)
\ (X\subseteq V(G)).
\end{equation}
In particular, $d_R(v)\ge d_G(v)/2$ for every $v\in V(G)$.
\end{lem}

\begin{proof}
We first characterize $\CS(G)^\perp$. The graph $G$ is connected by Proposition~\ref{pro3.3}, so we may fix a rooted spanning tree of $G$. Suppose that an edge set meets every cycle in an even number of edges, and color each vertex according to the parity of the number of selected edges on its tree path to the root. A tree edge is selected exactly when its endpoints have different colors. For a non-tree edge, the same assertion follows by applying the parity condition to its fundamental cycle. Thus the selected edges form a cut. Conversely, every cycle crosses a cut an even number of times. It follows that $\CS(G)^\perp$ consists precisely of the edge cuts.

Since $\HS(G)$ is a proper subspace of $\CS(G)$, the difference $\HS(G)^\perp\setminus\CS(G)^\perp$ is nonempty. Let $R$ be a spanning subgraph whose edge set has maximum cardinality in this difference. The symmetric difference with any edge cut remains in $\HS(G)^\perp\setminus\CS(G)^\perp$. For the cut between $X$ and $V(G)\setminus X$, this operation changes the cardinality by $e_G(X,V(G)\setminus X)-2e_R(X,V(G)\setminus X)$. By the maximality of $R$, this quantity is nonpositive, which proves~\eqref{eq-3}. To verify $R\ne G$, choose an edge of $G$. Proposition~\ref{pro3.3} provides a spanning path between its endpoints, and the path together with the chosen edge forms a Hamilton cycle. Its length is odd, so $E(G)$ is not orthogonal to this cycle. Hence $R\ne G$.
\end{proof}

For the next two lemmas, we first fix a subgraph $R$ as in Lemma~\ref{lem4.1}. The arguments adapt the short-path and cycle constructions of Hefetz and Krivelevich~\cite{hef}; both proofs are given below. We first use~\eqref{eq-3} to join two large balls in $R$ while avoiding prescribed vertices, and then use the resulting paths to construct an even cycle with exactly one edge outside $R$.

\begin{lem}\label{lem4.2}
Let $F\subseteq V(G)$ with $|F|\le100\hstar$, and let $x,y$ be distinct vertices of $V(G)\setminus(S\cup F)$. Then $R$ has an $x,y$-path of length at most $\hstar$ whose internal vertices lie in $U\setminus F$.
\end{lem}

\begin{proof}
Let $H=R[(U\setminus F)\cup\{x,y\}]$. Each vertex outside $S$ has at most one neighbor in $Z$ by~\ref{iteA2}, and $d_R(v)\ge d_G(v)/2$ by Lemma~\ref{lem4.1}. Together with~\ref{iteA1}, these bounds imply that each $v\in V(H)$ satisfies $d_H(v)\ge d_G(v)/2-1-|F|\ge L/25$ for sufficiently large $n$, since $\hstar=o(L)$. We may therefore apply Lemma~\ref{lem2.6} to $H$.

Consider the balls $B_H(x,h)$, for nonnegative integers $h$, as defined in Section~\ref{sec2}. By Lemma~\ref{lem2.6}, as long as a ball has size at most $\sstar$, its size increases by a factor of at least $r$ when its radius is increased by one. Thus a ball of radius at most $\lceil\log(\sstar+1)/\log r\rceil+1$ has more than $\sstar$ vertices. We apply Lemma~\ref{lem2.6} to any $\sstar$ vertices of this ball. The next ball then contains at least $a=r\sstar$ vertices, and its radius is at most $(2+o(1))L/\ell$. The parameters satisfy $a=o(n)$ and $a\gg\bstar$.

We next show that one further increase in the radius produces a ball with at least $0.48n$ vertices. Suppose that the enlarged ball has fewer vertices. Since $|V(H)|=n-o(n)$, there is a set $B_0\subseteq V(H)$ of size $\lceil0.51n\rceil$ outside it. Let $A_0$ be a set of $a$ vertices in the preceding ball. Any $R$-neighbor of a vertex of $A_0$ that lies in $V(H)$ belongs to the enlarged ball, so no edge of $R$ joins $A_0$ to $B_0$. It follows that
\begin{equation}\label{eq-4}
e_G(A_0,B_0)\le\sum_{v\in A_0}\bigl(d_G(v)-d_R(v)\bigr)
\le\frac12\sum_{v\in A_0}d_G(v).
\end{equation}
Both parts of any balanced partition of $A_0$ have size at least $\bstar$. An edge of $G[A_0]$ crosses a balanced partition chosen uniformly at random with probability $2\lfloor a/2\rfloor\lceil a/2\rceil/(a(a-1))$. The upper bound in~\ref{iteA5}, averaged over these partitions, therefore implies $2e_G(A_0)\le1.001p_0 a(a-1)$. We then apply~\ref{iteA5} to $A_0$ and its complement and obtain
\[
\sum_{v\in A_0}d_G(v)
=2e_G(A_0)+e_G(A_0,V(G)\setminus A_0)
\le1.001p_0 a(n-1).
\]
The lower bound in~\ref{iteA5} also yields $e_G(A_0,B_0)\ge0.999p_0 a\lceil0.51n\rceil$. Since $0.999\cdot0.51>1.001/2$, the last two estimates contradict~\eqref{eq-4}.

We repeat the argument around $y$. The two balls so obtained have size at least $0.48n$ and radius at most $(2+o(1))L/\ell+1$. If they intersect, paths from $x$ and $y$ to a common vertex form an $x,y$-walk of length less than $\hstar$. If they are disjoint, let $A,B$ be subsets of the respective balls, each of size $\lfloor0.4n\rfloor$, and let $C_0=V(G)\setminus(A\cup B)$. By~\eqref{eq-3} and~\ref{iteA5},
\begin{align*}
e_R(A,B)
&\ge\tfrac12e_G(A,V(G)\setminus A)-e_G(A,C_0)\\
&\ge\left(\tfrac12\cdot0.999\cdot0.4\cdot0.6
-1.001\cdot0.4\cdot0.2+o(1)\right)n^2p_0>0.
\end{align*}
We can therefore join the two balls by an edge of $R$ and obtain an $x,y$-walk of length at most $(4+o(1))L/\ell+3<\hstar$. In either case, a shortest $x,y$-path within the walk has all its internal vertices in $U\setminus F$, as required.
\end{proof}

\begin{lem}\label{lem4.3}
There is an even cycle $C$ in $G$ of length at most $4\hstar+3$ with exactly one edge outside $R$. Moreover, $C$ satisfies one of the following conditions.
\begin{enumerate}[label=(\roman*)]
\item $V(C)\cap S=\varnothing$ and $|V(C)\cap T|\le1$.
\item $V(C)\cap S=\{u\}$ and $V(C)\cap T=N_C(u)$.
\end{enumerate}
Here $N_C(u)$ denotes the two neighbors of $u$ on $C$.
\end{lem}

\begin{proof}
Recall from Definition~\ref{def2.1} that $T=N_G(S)$, $Z=S\cup T$ and $U=V(G)\setminus Z$. By Lemma~\ref{lem4.2}, $R[U]$ is connected and has diameter at most $\hstar$. For each $t\in T$, we have $d_R(t)\ge L/20$ and $d_G(t,Z)\le1$, so $t$ has an $R$-neighbor in $U$. Every vertex of $S$ has at least two $R$-neighbors, all of which lie in $T$. Thus $R$ is connected.

\emph{Case 1. The graph $R$ is bipartite.}
Let $(A',B')$ be the bipartition of $R$. All edges of $R$ join its two parts, but $E(R)$ is not a cut of $G$. We may therefore choose an edge $xy\in E(G)\setminus E(R)$ with one endpoint in each part. Suppose first that $x,y\notin S$. By Lemma~\ref{lem4.2}, we can join them by an $x,y$-path in $R$ whose internal vertices lie in $U$. This path has odd length because its endpoints lie in different parts. Together with $xy$, it forms an even cycle with exactly one edge outside $R$. The independence of $T$ in $G$ means that at most one of $x,y$ belongs to $T$, so the cycle satisfies (i).

Suppose next that one endpoint, say $x$, belongs to $S$. Its neighbor $y$ then lies in $T$. Since $x$ has at least two $R$-neighbors, choose $x'\in N_R(x)\setminus\{y\}$. The vertices $x'$ and $y$ belong to the same part of $R$, so by Lemma~\ref{lem4.2}, we can join them by an even path with all internal vertices in $U$. This path, together with $x'x$ and $xy$, forms an even cycle satisfying (ii). In both subcases the cycle has length at most $\hstar+2$.

\emph{Case 2. The graph $R[U]$ is nonbipartite.}
Choose a breadth-first spanning tree of $R[U]$, whose depth is at most $\hstar$. If every edge joined vertices at depths of opposite parity, these parities would define a proper two-coloring of $R[U]$. Since $R[U]$ is nonbipartite, some edge joins vertices whose depths have the same parity. This edge and the tree path between its endpoints form an odd cycle $Q$ of length at most $2\hstar+1$.

Let $xy\in E(G)\setminus E(R)$. Since $S$ is independent, we may name its endpoints so that $y\notin S$. If $x\in S$, choose $x'\in N_R(x)\setminus\{y\}$. Otherwise, let $x'=x$. We will join $x$ and $y$ to distinct vertices of $Q$ by vertex-disjoint paths that meet $Q$ only at their endpoints.

If $x'\in V(Q)$, let $u=x'$ and take the trivial path from $x'$ to $u$. Otherwise, choose $u\in V(Q)\setminus\{y\}$ and apply Lemma~\ref{lem4.2} to $x',u$, with forbidden set $(V(Q)\setminus\{u\})\cup\{y\}$. When $x\in S$, we prepend the edge $xx'$ to the resulting path. In either case we obtain an $x,u$-path that avoids $y$ and meets $Q$ only at $u$.

If $y\in V(Q)$, let $v=y$ and take the trivial path from $y$ to $v$. Otherwise, choose $v\in V(Q)\setminus\{u\}$ and apply Lemma~\ref{lem4.2} to $y,v$, forbidding all vertices of the first path and all vertices of $Q$ except $v$. The forbidden sets in these two applications have fewer than $100\hstar$ vertices and contain neither endpoint of the path to be constructed. Hence the resulting $y,v$-path is disjoint from the $x,u$-path and meets $Q$ only at $v$.

The two paths and the edge $xy$ now form a $u,v$-path of length at most $2\hstar+2$, with exactly one edge outside $R$ and no internal vertex in $Q$. The two $u,v$-arcs of the odd cycle $Q$ have opposite parities, so one of them completes this path to an even cycle of length at most $4\hstar+3$. If $x\notin S$, its only possible vertices in $T$ are $x,y$, and the independence of $T$ implies that at most one of them belongs to $T$. If $x\in S$, its only vertices in $Z$ are $x,x',y$, with $x',y$ as the two cycle neighbors of $x$. These are precisely conditions (i) and (ii), respectively.

\emph{Case 3. The graph $R[U]$ is bipartite and $R$ is nonbipartite.}
Let $(A,B)$ be the bipartition of $R[U]$. We first construct an even path in $R$ from some $a\in A$ to some $b\in B$. Its internal vertices will consist of one vertex of $T$, or of one vertex of $S$ together with two of its neighbors in $T$.

If some $t\in T$ has $R$-neighbors $a\in A$ and $b\in B$, we use the path $a,t,b$. Otherwise, the $R$-neighbors in $U$ of each vertex of $T$ all lie in one part. We extend the proper two-coloring of $R[U]$ to $T$ by assigning each such vertex the color opposite to that of its $R$-neighbors in $U$. The independence of $T$ ensures that the resulting coloring of $R[U\cup T]$ is proper. Since $S$ is independent and $R$ is nonbipartite, some $u\in S$ must have two $R$-neighbors $t,t'\in T$ of different colors. Let $a\in N_R(t)\cap U$ and $b\in N_R(t')\cap U$. These vertices also have different colors, and we may name the parts so that $a\in A$ and $b\in B$. We then use the path $a,t,u,t',b$.

We may interchange the two parts and reverse the path, if necessary, to assume that $|A|\ge|B|$. Since $|U|\ge n-\sqrt n$ and $\bstar=o(n)$, there is a partition of $A\setminus\{a\}$ into two sets of size at least $\bstar$. By~\ref{iteA5}, there is an edge $vw$ between these sets. Both endpoints lie in the same part of $R[U]$, so $vw$ lies outside $R$. We apply Lemma~\ref{lem4.2} first to obtain an $a,v$-path in $R[U]$ that avoids $b,w$, and then to obtain a $b,w$-path in $R[U]$ that avoids all vertices of the first path. Each path has length at most $\hstar$. The $a,v$-path has even length, whereas the $b,w$-path has odd length. Together with $vw$ and the even path through $Z$, they form an even cycle of length at most $2\hstar+5\le4\hstar+3$. Its only edge outside $R$ is $vw$. It satisfies (i) when the even path through $Z$ has its internal vertex in $T$, and (ii) when this path contains a vertex of $S$.
\end{proof}

We use the parity-switcher construction of Christoph, Nenadov and Petrova~\cite{chr} to obtain two paths with the same endpoints and vertex set whose edge sets have symmetric difference equal to the cycle from Lemma~\ref{lem4.3}. The restrictions on that cycle ensure that one common complementary path completes both of them to Hamilton cycles.

\begin{prop}\label{pro4.4}
The graph $G$ satisfies $\HS(G)=\CS(G)$.
\end{prop}

\begin{proof}
Suppose for a contradiction that $\HS(G)\ne\CS(G)$, and choose $R$ and $C$ as in Lemmas~\ref{lem4.1} and~\ref{lem4.3}, respectively. Let $2k$ be the length of $C$. Since $C$ contains at most one vertex of $S$, we may choose two opposite vertices $x,y$ of $C$ outside $S$. Label the vertices of the two arcs from $x$ to $y$ so that the cyclic order on $C$ is $x,a_1,\ldots,a_{k-1},y,b_{k-1},\ldots,b_1$. We also have $2k\le4\hstar+3=O(L/\ell)$.

\textbf{Step 1. Preserve neighbors of the exceptional vertices.}
For each $v\in V(C)\setminus S$, let $v'=v$. If $C$ contains a vertex $u\in S$, the only vertices of $T$ on $C$ are its two cycle neighbors. Since $d_G(u)\ge3$, we may choose a representative $u'\in N_G(u)\setminus V(C)$. The representatives $v'$ of the vertices of $C$ are therefore all distinct. Let
\[
\mathcal T_0=\{v':v\in V(C)\},\ W_0=V(C)\cup \mathcal T_0,\ Y_0=V(G)\setminus(Z\cup W_0).
\]
Each vertex outside $S$ has at most one neighbor in $Z$, and $|W_0|=O(L/\ell)$. Hence, for sufficiently large $n$, every $v\notin S$ satisfies $d_G(v,Y_0)\ge L/10-1-|W_0|\ge3L/40$. We apply Lemma~\ref{lem3.1} with $Q=V(G)\setminus S$ and $\gamma=3/40$ and obtain a partition $Y_0=A_0\mathbin{\dot\cup}B_0$ such that every vertex outside $S$ has at least $L/40$ neighbors in each part. Since $|Q|\ge n/2$, the size bound in that lemma also shows that each part has size at least $n/800$.

Let $A=A_0\cup(Z\setminus W_0)$ and $B=B_0\cup \mathcal T_0$. These sets are disjoint. Every vertex of $B$ lies outside $S$, so it has at least $L/40$ neighbors in $B_0$.

\textbf{Step 2. Construct the two paths.}
Let $G'=G[B\setminus\{x,y\}]-E(C)$. Each vertex loses at most two neighbors when we delete $x,y$ and at most two incident edges when we delete $E(C)$. The minimum degree of $G'$ is therefore at least $L/40-4\ge L/50$ for sufficiently large $n$. Let $\mathcal T=\mathcal T_0\setminus\{x,y\}$. Since $|\mathcal T|=2k-2=O(L/\ell)=o(L)$, we have $d_{G'}(v,\mathcal T)\le d_{G'}(v)/10$ for all $v\in V(G')$.

The number of pairs is $k-1\le\sqrt n$ for sufficiently large $n$, so by Lemma~\ref{lem2.8}, we can find pairwise vertex-disjoint paths in $G'$ joining $a_i'$ to $b_i'$ for $1\le i\le k-1$. If $u\in S\cap V(C)$, we extend the path ending at $u'$ by the edge $u'u$. Denote the resulting paths by $P_i$, with endpoints $a_i$ and $b_i$ for $1\le i\le k-1$. They remain pairwise vertex-disjoint. Each ordinary vertex on $C$ other than $x,y$ is a prescribed endpoint, so it cannot be an internal vertex of any of these paths. The vertices $x,y$ are absent from $G'$, and the only possible exceptional vertex on $C$ is added at its own endpoint. Thus the paths have no internal vertices on $C$. They also avoid $E(C)$, since those edges were deleted from $G'$ and the possible extension joins $u$ to a vertex $u'$ outside $C$.

Let $H_C$ be the union of $C$ and the paths $P_i$, and let $M_0,M_1$ be the two alternating perfect matchings of $C$. For $j\in\{0,1\}$, let $Q_j$ be the union of $M_j$ and all the paths $P_i$. We verify that each $Q_j$ is a spanning $x,y$-path of $H_C$. After we contract each $P_i$ to a vertex $w_i$, both arcs of $C$ from $x$ to $y$ follow the same vertex sequence $x,w_1,\ldots,w_{k-1},y$. The first edges of the two arcs belong to different matchings, and the edges on each arc alternate between $M_0$ and $M_1$. Thus each matching supplies exactly one edge between each pair of consecutive vertices in this sequence. It follows that $Q_j$ is connected. In the original graph $Q_j$, the vertices $x,y$ have degree one and every other vertex has degree two. Hence $Q_j$ is a spanning $x,y$-path of $H_C$. Since the paths $P_i$ occur in both $Q_0$ and $Q_1$, whereas the two matchings partition $E(C)$, we also have
\begin{equation}\label{eq-5}
E(Q_0)\mathbin{\triangle}E(Q_1)=E(C).
\end{equation}

Figure~\ref{fig3} illustrates the construction for $k=3$. Solid segments are edges of $C$, dashed paths are $P_1,P_2$, and the dotted path $P$ is the common complementary path constructed in the next step. The heavier edge $e$ is the unique edge of $C$ outside $R$. The common paths cancel in the symmetric difference of the two Hamilton cycles, leaving $E(C)$.

\begin{center}
\begin{minipage}{\linewidth}
\centering
\begin{tikzpicture}[
  x=1cm,y=1cm,
  every node/.style={font=\footnotesize},
  vertex/.style={circle,draw=black,fill=white,inner sep=0pt,minimum size=6mm,line width=0.55pt},
  point/.style={circle,draw=black,fill=white,inner sep=0pt,minimum size=3.5pt,line width=0.55pt},
  cycle edge/.style={line width=0.7pt},
  marked edge/.style={line width=1.4pt},
  connecting path/.style={line width=0.8pt,dash pattern=on 3pt off 2pt},
  complementary path/.style={line width=0.85pt,dash pattern=on 0.6pt off 1.8pt}
]
\foreach \panel/\shift in {0/0,1/5.1,2/10.2}{
  \begin{scope}[xshift=\shift cm]
  \coordinate (x) at (0,0);
  \coordinate (a) at (0.85,0.95);
  \coordinate (b) at (2.55,0.95);
  \coordinate (y) at (3.4,0);
  \coordinate (c) at (2.55,-0.95);
  \coordinate (d) at (0.85,-0.95);
  \coordinate (pa) at (0.85,0);
  \coordinate (pb) at (2.55,0);
  \draw[connecting path] (a)--(d);
  \draw[connecting path] (b)--(c);
  \node[fill=white,inner sep=1pt] at (0.52,0) {$P_1$};
  \node[fill=white,inner sep=1pt] at (2.88,0) {$P_2$};
  \ifnum\panel=0
    \draw[cycle edge] (x)--(a)--(b)--(y)--(c)--(d)--cycle;
    \draw[marked edge] (x)--(a);
    \node at (0.10,0.72) {$e$};
  \fi
  \ifnum\panel=1
    \draw[cycle edge] (x)--(a) (b)--(y) (c)--(d);
    \draw[marked edge] (x)--(a);
    \node at (0.10,0.72) {$e$};
  \fi
  \ifnum\panel=2
    \draw[cycle edge] (a)--(b) (y)--(c) (d)--(x);
  \fi
  \ifnum\panel>0
    \draw[complementary path] (x) .. controls (0.05,-2.3) and (3.35,-2.3) .. (y);
    \node[point] at (1.7,-1.725) {};
    \node[fill=white,inner sep=1pt] at (1.7,-2.02) {$P$};
  \fi
  \foreach \v/\lab in {x/x,a/a_1,b/a_2,y/y,c/b_2,d/b_1}{
    \node[vertex] at (\v) {$\lab$};
  }
  \foreach \v in {pa,pb}{
    \node[point] at (\v) {};
  }
  \ifnum\panel=0
    \node at (1.7,-2.45) {(a) $C\cup P_1\cup P_2$};
  \fi
  \ifnum\panel=1
    \node at (1.7,-2.45) {(b) $P\cup Q_0$};
  \fi
  \ifnum\panel=2
    \node at (1.7,-2.45) {(c) $P\cup Q_1$};
  \fi
  \end{scope}
}
\end{tikzpicture}
\captionsetup{font=small,labelsep=period,hypcap=false}

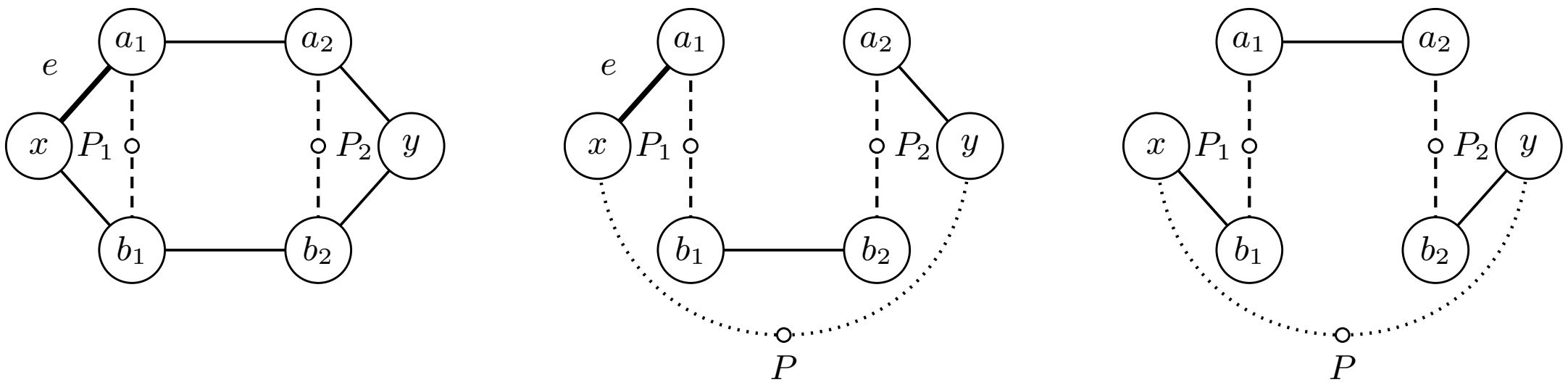
\captionof{figure}{Two Hamilton cycles obtained from a parity-switcher with $k=3$.}
\label{fig3}
\end{minipage}
\end{center}

\textbf{Step 3. Complete both paths.}
Let $W=V(G)\setminus(V(H_C)\setminus\{x,y\})$. The vertices of $C$ lie in $W_0$, and the internal vertices of the paths $P_i$ lie in $B$. Since $A$ is disjoint from both $W_0$ and $B$, we have $A\subseteq W$. In particular, $|W|\ge n/800$. Every vertex of $W\setminus S$ still has at least $L/40$ neighbors in $A_0$.

We check the remaining hypothesis of Lemma~\ref{lem3.2} for vertices in $W\cap S$. If $V(C)\cap S=\varnothing$, then $W_0=V(C)$ and $|W_0\cap T|\le1$. Each vertex of $S$ therefore has at most one neighbor in $W_0$. If $V(C)\cap S=\{u\}$, every vertex of $W_0\cap T$ is a neighbor of $u$, so~\ref{iteA2} excludes adjacency to any vertex of $S\setminus\{u\}$. Thus, in both cases, each vertex of $W\cap S$ has at least two neighbors in $Z\setminus W_0\subseteq A$. These neighbors are distinct from $x,y\in W_0$. We may now apply Lemma~\ref{lem3.2} to obtain a spanning $x,y$-path $P$ of $G[W]$.

For each $j\in\{0,1\}$, the paths $P$ and $Q_j$ meet only at $x,y$ and together cover $V(G)$. Their union is therefore a Hamilton cycle of $G$. The common path $P$ cancels in the symmetric difference of these two cycles, leaving $E(C)$ by~\eqref{eq-5}. Since $C$ has even length and exactly one edge outside $R$, the number $|E(C)\cap E(R)|$ is odd. Consequently, the numbers of edges in $R$ on the two Hamilton cycles have opposite parities. One of them has an odd number of edges in $R$, contrary to Lemma~\ref{lem4.1}.
\end{proof}

Propositions~\ref{pro3.3} and~\ref{pro4.4} complete the proof of Theorem~\ref{the2.2}.

\section{Verification at the hitting time}
\label{sec5}

The edge count at $\tau_3$ is random. We therefore couple the process with two binomial random graphs whose edge counts bound $\tau_3$ from below and above, and verify the hypotheses of Theorem~\ref{the2.2} for every intermediate graph of minimum degree at least three.

Recall that $[n]=\{1,\ldots,n\}$. We assign independent random labels $\xi_e$, each uniformly distributed on $[0,1]$, to the edges $e\in E(K_n)$. Almost surely, these labels are distinct and order the edges in a uniformly random permutation. The first $m$ edges then form $G(n,m)$. For $p\in[0,1]$, let $\mathbf G(p)$ be the graph on $[n]$ with edge set $\{e\in E(K_n):\xi_e\le p\}$. Each edge belongs to this graph independently with probability $p$, so $\mathbf G(p)$ has distribution $\mathbb{G}(n,p)$.

Let $\omega=\log\ell$ and $p_\pm=(L+2\ell\pm\omega)/n$, and write $G^-=\mathbf G(p_-)$ and $G^+=\mathbf G(p_+)$. For $m_\pm=|E(G^\pm)|$, this coupling gives $G^-=G(n,m_-)$ and $G^+=G(n,m_+)$. We have $p_\pm/p_0\to1$, where $p_0$ is given in~\eqref{eq-1}. Throughout this section, let
\begin{equation}\label{eq-6}
S=\{v:d_{G^-}(v)\le L/10\}.
\end{equation}

\begin{lem}\label{lem5.1}
With probability $1-o(1)$, every spanning graph $G$ with $G^-\subseteq G\subseteq G^+$ and $\delta(G)\ge3$ satisfies the criterion in Definition~\ref{def2.1} with the set $S$ in~\eqref{eq-6}.
\end{lem}

\begin{proof}
An upper bound on a degree or an edge count in $G^+$ holds in every spanning subgraph, while a lower bound in $G^-$ holds in every spanning supergraph. It therefore suffices to prove the corresponding estimates for these two graphs.

\emph{Step 1. Degrees and separation.}
For sufficiently large $n$, we have $p_+\le2L/n$. The second inequality in~\eqref{eq-chernoff-one-sided} gives $\Prb(d_{G^+}(v)\ge10L)\le(e/5)^{10L}=o(n^{-1})$, and the union bound yields $\Delta(G^+)\le10L$ with high probability. The definition of $S$ also ensures that each vertex outside $S$ has degree greater than $L/10$ in $G^-$.

Let $I=0.9+0.1\log0.1>3/5$. The mean degree in $G^-$ is $L+o(L)$. The first inequality in~\eqref{eq-chernoff-one-sided}, applied to the probability of a degree at most $L/10$, gives
\begin{equation}\label{eq-7}
\Prb(v\in S)\le n^{-I+o(1)}.
\end{equation}
It follows that $\Ex|S|\le n^{1-I+o(1)}$. By Markov's inequality, $|S|\le n^{0.45}$ with high probability. Together with $\Delta(G^+)\le10L$, this gives $|S\cup N_{G^+}(S)|\le(1+10L)n^{0.45}\le\sqrt n$.

The separation conditions follow by excluding short paths between vertices of $S$ and short cycles through $S$. Fix a simple path of length $j\in\{1,2,3,4\}$ in $K_n$. For both endpoints to belong to $S$, each must have at most $L/10$ neighbors in $G^-$ outside the path. These two counts depend on disjoint sets of edge labels, which are also disjoint from the labels on the path. Each count excludes only a bounded number of possible neighbors, so the estimate~\eqref{eq-7} for a degree at most $L/10$ still applies. The union bound over all such paths is therefore
\[
\sum_{j=1}^4 n^{j+1}p_+^j n^{-2I+o(1)}=n^{1-2I+o(1)}=o(1).
\]
For a triangle or a four-cycle through $S$, we instead specify one vertex on the cycle that belongs to $S$. The same argument gives the bound
\[
\sum_{j=3}^4 jn^jp_+^j n^{-I+o(1)}=o(1)
\]
on the probability that a triangle or a four-cycle of $G^+$ contains a vertex of $S$.

Suppose that $|S\cup N_{G^+}(S)|\le\sqrt n$ and that none of these path or cycle configurations occurs. We check~\ref{iteA2} for any spanning graph $G$ with $G^-\subseteq G\subseteq G^+$. The size bound gives $|S\cup N_G(S)|\le\sqrt n$. Since no path of length at most four joins distinct vertices of $S$, the set $S$ is independent. If two vertices of $T=N_G(S)$ were adjacent, their neighbors in $S$ would give either a path of length three between distinct vertices of $S$ or a triangle containing a vertex of $S$. It follows that $T$ is independent as well. To prove the remaining degree bound, suppose that $v\notin S$ has two distinct neighbors in $Z=S\cup T$. Associate each of these neighbors in $T$ with one of its neighbors in $S$ and each of them in $S$ with itself. If the associated vertices are distinct, there is a walk between them of length at most four, and hence a path of at most that length. If the associated vertices coincide, the two neighbors of $v$ give a four-cycle through the common vertex of $S$ when both lie in $T$, and a triangle when one lies in $S$. Each possibility contradicts the assumed absence of these configurations.

\emph{Step 2. Small sets.}
For $1\le a\le\sstar$, write $h=\log(n/a)$. Then $h\ge\ell-2\log\ell$. We apply the binomial bound $\Prb(\operatorname{Bin}(M,p)\ge t)\le(eMp/t)^t$ with $t=\lceil a\qstar\rceil$. When $t>M$, the event has probability zero. In the other case, the bases in the estimates below are less than one for sufficiently large $n$. We may therefore replace the exponent $t$ by $a\qstar$ in each upper bound.

The union bound over vertex sets of size $a$ gives the following estimate for~\ref{iteA3}.
\begin{align*}
\Prb\bigl(\exists X\text{ with }|X|=a\text{ and }e_{G^+}(X)>a\qstar\bigr)
&\le\left(\frac{en}{a}\right)^a
      \left(\frac{ea\ell}{n}\right)^{a\qstar}\\
&=\exp\bigl(a[1+h-\qstar(h-\log\ell-1)]\bigr)\\
&\le\exp(-a\qstar h/2),
\end{align*}
where the last inequality holds uniformly for $1\le a\le\sstar$.

For~\ref{iteA4}, it suffices to consider $|Y|=ra$. To justify this reduction, note that $(r+1)\sstar=o(n)$, so any set disjoint from $X$ of size less than $ra$ can be enlarged to size $ra$ while remaining disjoint from $X$. The union bound over disjoint pairs gives
\begin{align*}
&\Prb\bigl(\exists X,Y\text{ disjoint with }|X|=a,\ |Y|=ra
   \text{ and }e_{G^+}(X,Y)>a\qstar\bigr)\\
&\le
\left(\frac{en}{a}\right)^a
\left(\frac{en}{ra}\right)^{ra}
\left(\frac{2era\ell}{n}\right)^{a\qstar}\\
&=
\exp\bigl(a[1+h+r(1+h-\log r)
-\qstar(h-\log r-\log\ell-\log(2e))]\bigr)\\
&\le\exp(-0.3a\qstar h).
\end{align*}
The last inequality uses $(\log r+\log\ell+\log(2e))/h\le1/2+o(1)$ uniformly over the given range. Also, $r/\qstar=o(1)$ gives $1+h+r(1+h-\log r)=o(\qstar h)$. Since $\qstar h\ge L/2$ for sufficiently large $n$, each of the two bounds sums to $o(1)$ over $a$.

\emph{Step 3. Large sets.}
Let $\varepsilon=10^{-4}$, and fix disjoint sets $X,Y$ of sizes $a,c\ge\bstar$, respectively. The edge count between them is binomial in each of $G^-$ and $G^+$, so~\eqref{eq-chernoff} gives
\[
\Prb\bigl(|e_{G^\pm}(X,Y)-p_\pm ac|>\varepsilon p_\pm ac\bigr)\le2\exp(-c_1p_0 ac)
\]
for some constant $c_1>0$.

For fixed $a,c$, there are at most $\binom na\binom nc\le\exp(2(a+c)\ell)$ choices of $X,Y$. The estimate $p_0 ac\ge\tfrac12(a+c)p_0\bstar\ge\tfrac12(a+c)\ell^{3/2}$ shows that the preceding tail bound absorbs this factor. For sufficiently large $n$, the probability of a violation for some pair of these sizes in either graph is therefore at most $4\exp(-c_1p_0 ac/2)$. Furthermore, $p_0 ac\ge p_0\bstar^2$, and $p_0\bstar^2$ is of order $n\ell^3/L\gg\log n$. The union bound over the at most $n^2$ choices of $a,c$ is thus $o(1)$. Since $p_\pm/p_0\to1$, the upper estimate in $G^+$ and the lower estimate in $G^-$ imply~\eqref{eq-2} for every spanning graph between $G^-$ and $G^+$.

On the intersection of the events proved above, conditions~\ref{iteA2}--\ref{iteA5} and the degree bounds in~\ref{iteA1} all hold. The remaining part of~\ref{iteA1} is the hypothesis $\delta(G)\ge3$.
\end{proof}

\begin{lem}\label{lem5.2}
With probability $1-o(1)$, we have $\delta(G^-)\le2$ and $\delta(G^+)\ge3$.
\end{lem}

\begin{proof}
For the lower graph, we count the vertices whose degree is exactly two. Let $X_2$ denote their number in $G^-$ and let $\mu=\Ex X_2$. The binomial formula gives
\[
\mu=n\binom{n-1}{2}p_-^2(1-p_-)^{n-3}\sim\tfrac12e^\omega=\tfrac12\ell\longrightarrow\infty.
\]
Let $a_i=\Prb(\operatorname{Bin}(n-2,p_-)=i)$ for $i=1,2$. For two distinct vertices, the degree counts that exclude the edge between them are independent. Conditional on the absence or presence of this edge, a given vertex has degree two with probability $a_2$ or $a_1$, respectively. Its unconditional probability is therefore $\pi=(1-p_-)a_2+p_-a_1$, and the covariance of the two degree-two indicators is
\[
(1-p_-)a_2^2+p_-a_1^2-\pi^2=p_-(1-p_-)(a_2-a_1)^2.
\]
The ratio $a_1/a_2=2(1-p_-)/((n-3)p_-)=O(1/L)$ bounds this covariance by $2p_-\pi^2$ for sufficiently large $n$. Hence $\Var X_2\le\mu+2p_-\mu^2$. Chebyshev's inequality then gives $\Prb(X_2=0)\le\mu^{-1}+2p_-=o(1)$.

For the upper graph $G^+$, a first-moment bound suffices. The expected number of vertices of degree at most two is
\[
n\sum_{j=0}^2\binom{n-1}{j}p_+^j(1-p_+)^{n-1-j}\sim\tfrac12e^{-\omega}=\frac1{2\ell}=o(1).
\]
By Markov's inequality, $\delta(G^+)\ge3$ with high probability.
\end{proof}

We use the following consequence of a monotonicity result of Heinig~\cite{hei1} to extend the conclusion from $\tau_3$ to every subsequent step. We include the short proof.

\begin{lem}\label{lem5.3}
If $G$ is Hamilton-connected and satisfies $\HS(G)=\CS(G)$, then every spanning supergraph $G'$ of $G$ is Hamilton-connected and satisfies $\HS(G')=\CS(G')$.
\end{lem}

\begin{proof}
It suffices to consider the addition of one missing edge $xy$. All spanning paths of $G$ remain, so the resulting graph is Hamilton-connected. Let $P$ be a spanning $x,y$-path in $G$. The new edge completes it to a Hamilton cycle $Q=P+xy$. For any cycle $C$ containing $xy$, this edge cancels in $E(C)\mathbin{\triangle}E(Q)$, leaving an even-degree subgraph of $G$. Its edge set belongs to $\CS(G)$ and is therefore a sum of Hamilton cycles of $G$. The addition of $E(Q)$ to this sum expresses $E(C)$ as a sum of Hamilton cycles in the resulting graph. Every cycle that avoids $xy$ already belongs to $\CS(G)$.
\end{proof}

\begin{proof}[\textbf{Proof of Theorem~\ref{the1.1}}]
By Lemmas~\ref{lem5.1} and~\ref{lem5.2}, all their conclusions hold simultaneously with probability $1-o(1)$. On this event, the identities $G^-=G(n,m_-)$ and $G^+=G(n,m_+)$, together with Lemma~\ref{lem5.2}, give $m_-<\tau_3\le m_+$.

Thus $G^-\subseteq G(n,\tau_3)\subseteq G^+$ and $\delta(G(n,\tau_3))\ge3$. By Lemma~\ref{lem5.1}, the graph at $\tau_3$ satisfies the deterministic criterion. By Theorem~\ref{the2.2}, $G(n,\tau_3)$ is therefore Hamilton-connected and satisfies $\HS(G(n,\tau_3))=\CS(G(n,\tau_3))$. Each $G(n,m)$ with $m\ge\tau_3$ is a spanning supergraph of $G(n,\tau_3)$, and by Lemma~\ref{lem5.3}, both properties hold for all such $m$.
\end{proof}
\section{Declaration of competing interest}
The authors declare that they have no known competing financial interests or personal relationships that could have appeared to influence the work reported in this paper.
\section{Data availability}
No data were used for the research described in this article.
\section{AI Declaration}
During the preparation of this manuscript, the authors used AI tools to assist in developing mathematical arguments, checking calculations, and improving the exposition. The authors directed this process, independently checked the proofs, and take full responsibility for the mathematical content and presentation.


\begin{thebibliography}{99}
\small

\bibitem{ajt} M. Ajtai, J. Koml\'os, E. Szemer\'edi, First occurrence of Hamilton cycles in random graphs, in: Cycles in Graphs, North-Holland Math. Stud. 115, North-Holland, Amsterdam, 1985, pp. 173--178.

\bibitem{alo} N. Alon, J.H. Spencer, The Probabilistic Method, 4th ed., John Wiley \& Sons, Hoboken, NJ, 2016.

\bibitem{als} B. Alspach, S.C. Locke, D. Witte, The Hamilton spaces of Cayley graphs on abelian groups, Discrete Math. 82 (1990) 113--126.

\bibitem{bar} J.D. Baron, J. Kahn, On the cycle space of a random graph, Random Structures Algorithms 54 (2019) 39--68.

\bibitem{bol1} B. Bollob\'as, The evolution of sparse graphs, in: Graph Theory and Combinatorics, Academic Press, London, 1984, pp. 35--57.

\bibitem{bol2} B. Bollob\'as, A. Thomason, Random graphs of small order, in: Random Graphs '83, North-Holland Math. Stud. 118, North-Holland, Amsterdam, 1985, pp. 47--97.

\bibitem{chr} M. Christoph, R. Nenadov, K. Petrova, The Hamilton space of pseudorandom graphs, J. Combin. Theory Ser. B 176 (2026) 254--267.

\bibitem{dem} B. DeMarco, A. Hamm, J. Kahn, On the triangle space of a random graph, J. Comb. 4 (2013) 229--249.

\bibitem{dra1} N. Dragani\'c, M. Krivelevich, R. Nenadov, Rolling backwards can move you forward: on embedding problems in sparse expanders, Trans. Amer. Math. Soc. 375 (2022) 5195--5216.

\bibitem{dra2} N. Dragani\'c, R. Montgomery, D. Munh\'{a} Correia, A. Pokrovskiy, B. Sudakov, Hamiltonicity of expanders: optimal bounds and applications, arXiv:2402.06603v2 (2024).

\bibitem{erd1} P. Erd\H{o}s, A. R\'enyi, On random graphs. I, Publ. Math. Debrecen 6 (1959) 290--297.

\bibitem{erd2} P. Erd\H{o}s, A. R\'enyi, On the evolution of random graphs, Publ. Math. Inst. Hungar. Acad. Sci. 5 (1960) 17--61.

\bibitem{gil} E.N. Gilbert, Random graphs, Ann. Math. Statist. 30 (1959) 1141--1144.

\bibitem{gle} R. Glebov, M. Krivelevich, On the number of Hamilton cycles in sparse random graphs, SIAM J. Discrete Math. 27 (2013) 27--42.

\bibitem{hef} D. Hefetz, M. Krivelevich, The Hamilton cycle space of random graphs, arXiv:2506.19731v2 (2025).

\bibitem{hei1} P.C. Heinig, On prisms, M\"obius ladders and the cycle space of dense graphs, European J. Combin. 36 (2014) 503--530.

\bibitem{hei2} P.C. Heinig, When Hamilton circuits generate the cycle space of a random graph, arXiv:1303.0026 (2013).



\bibitem{mon} R. Montgomery, Hamiltonicity in random graphs is born resilient, J. Combin. Theory Ser. B 139 (2019) 316--341.

\end{thebibliography}
\end{document}